ORESTA

# OPERATIONAL SCHEDULING OF OIL PRODUCTS PIPELINE WITH INTERMEDIATE EVENT OCCURRENCES

Seyyed Hamed Moghimi [1]*, Jafar Habibi [1], Hamid Jahad [1], Mohammad Amin Fazli[1]

[1] Department of Computer Engineering, Sharif University of Technology, Tehran, Iran

**Abstract:** *Oil products are the main source of energy in the world today. Distribution of these products is one of the main issues in the industry. The main tools for this work are pipelines, and along with it, railways, shipping and roads are also used. Optimal planning of pipelines is an example of decision-making problem and was the focus of many researchers in the past years. The use of mixed integer linear programming (MILP) is one of the efficient methods to solve this problem. However, models still ignore important operational challenges. Vulnerability to deal with incidents as well as lack of attention to other transportation methods as a complement to the pipeline are among the weak points of the existing models. In this research, we intend to facilitate the decision-making process for experts in the field of distribution of oil products. For this purpose, we must improve the existing MILP methods and modify them for use in the real operational environment in such a way that sufficient flexibility, the possibility of responding to incidents, and the ability to revise the program are added to them.*

## 1. Introduction

In the last two centuries, one of the main issues in the world has been the provision of energy. Oil and its refined products are considered as the most important fuel sources in the world. One on the important problems about oil products is their distribution. Considering the mechanical and chemical properties of these materials, such as their fluidity, density, as well as their flammability, their transportation is more sensitive, and for this reason, the related costs are also higher. Main tools used to move oil products today are pipelines, vehicles, railways, and vessels. Among those, pipelines have efficiency, higher safety and lower price and that is why many pipelines have been built around the world. For example, nearly two-thirds of oil products in the United States are distributed through pipelines.

* Corresponding author.
hmoghimi@ce.sharif.edu (S.H. Moghimi), jhabibi@sharif.edu (J. Habibi),
hamid.jahad@ce.sharif.edu (H. Jahad), fazli@sharif.edu (M.A. Fazli)

Optimal planning in order to make maximum use of the transmission capacity of pipelines is one of the main problems in this industry, which has allocated a significant amount of past and present researches. One of the most important points about this optimization problem is the very high financial turnover that exists in this regard. Considering the global price of oil products, the costs that may be imposed on a system in case of non-optimal planning of pipelines can be about several million dollars per day. Another reason is the strategic importance of fuel distribution. In all the countries of the world, energy carriers are considered as strategic assets, and providing the required amounts for normal use and saving for critical times is considered as one of the strategic goals of the administrations.

Answering the optimization problem is not an obvious and simple task, because the distribution network of oil products includes several pipes, different types of tanks and depots, refineries, pump stations and control centers, which together form a complex graph of different nodes and edges with unique characteristics and rules. Furthermore, various oil products are flowing inside the pipes, with different volumes based on the needs of destinations. Several rules, including physical constraints and safety principles, are also effective in making decisions about the order and amount of injection of different products. Therefore, it can be concluded that the nature of such a problem will not be simple and planning for it is a difficult and challenging activity.

The main tool used for solving the optimization problem of oil products distribution is mixed-integer linear programming (MILP) models. This method is developed in the last two decades and several researches have focused to solve the problems in this area. In this research, we proposed a novel method to take environmental changes and event occurrences into consideration. Our model can adapt to new environment and continues to execute the distribution plan, event different types of incidents occur and change the surrounding parameters.

## 2. Related Works

### 2.1. Early Researches

As of the beginning of the 21st century, researches were carried out in line with the use of the linear programming model for the optimization of multi-product pipelines. The first model presented for this purpose looked at the problem from a discrete point of view (Rejowski & Pinto, 2003). This means that the planning time horizon was divided into small intervals and also, the length of the pipelines were categorized into small volume intervals. This work made it easier to look at the problem in a discrete way, to understand mathematical relationships and equations and inequalities, but on the other hand, both the number of relations and computational complexity of the model increased significantly, and the accuracy and quality of the obtained results were not at a good level. Following that, other activities such as (Neiro & Pinto, 2004) and (Rejowski & Pinto, 2004) presented linear programming models with the same discrete approach.

With a little distance from these researches, the first linear programming model with continuous view was presented (Cafaro & Cerda, 2004). Continuous view versus discrete view means that instead of breaking time and place into small intervals and building model variables from the values of these intervals, the time and place of batches at key frames are considered as model variables. Key frames mean the moments when changes are to be made in the operation environment, for example a new batch is pumped into the pipeline or a batch is received at a station. With this point of view, the number of variables and constraints is significantly reduced and, as a result, the computational performance of the model grows up. The reason for reducing the dimensions of the problem in this model is that in discrete modeling, there are many times when the environment is operating normally in several consecutive time intervals, and practically all these intervals can be assumed as one key frame. After these early researches, a growing trend of researches began with the aim of using linear programming in modeling pipelines for more complex situations and involving more parameters into the problem to be closer to reality.

### 2.2. Topologies

In (Rejowski & Pinto, 2003) and (Cafaro & Cerda, 2004), both approaches consider only one refinery as the pipeline input and one or more stations as product depots. After that some researches focus on the variation of pipeline topology. In (MirHassani & Ghorbanalizadeh, 2008), (Mostafaei et al., 2015), and (Taherkhani et al., 2017) tree shaped pipelines are modeled. These tree topologies include the main pipeline and another level of child branches. More complicated topologies have been modeled in (Cafaro & Cerda, 2012), which introduced pipelines with mesh structure.

### 2.3. Input and output stations

One of the parameters considered in researches is number of input and output in the pipeline. Early models, consider the refinery as an entry point to the pipeline and some consecutive depots along the pipeline. This approach is also used in many later researches such as (Rejowski & Pinto, 2004), (Rejowski & Pinto, 2008), (Cafaro & Cerda, 2008), (MirHassani & Ghorbanalizadeh, 2008), (MirHassani & Fani Jahromi, 2011).

In the other hand, some researches simplify this assumption and just model one input and only one output terminal. This change results in less complexity and therefore, the possibility to other parameters being modeled (Relvas et al., 2006), (Moradi & MirHassani, 2014), and (Moradi & MirHassani, 2016).

More complex input and output terminals are also modeled. For examples (Cerda & Cafaro, 2009), (Cafaro & Cerda, 2010), (Magatao et al., 2015), (Mostafaei et al., 2015), (Mostafaei et al., 2016), and (Taherkhani, 2020) take into account multiple input nodes and dual-purposed terminals which can both pump new batches and withdraw product from them.

### 2.4. Scheduling time horizon

In all the researches around pipeline scheduling, one of the important parameters is the time horizon of planning, indicating the duration of the entire period for which planning is done. The importance of this parameter's value is related to a trade-off between the computational complexity of the model and its application in the real world. The shorter the time horizon results in fewer parameters and formulas and S.H.Moghimi et al. /Oper. Res. Eng. Sci. Theor. Appl. 3 (x) (2020) 1-20therefore, a less complex model. On the other hand, a shorter time horizon causes the future state of pipeline not to be considered. So it is possible that the result of the execution of models with shorter time horizon causes the pipeline network not to perform well in the long term, although they have good results in short time.

Researches have considered the time horizon as a determining parameter in research classification. There is a focus on three types of time horizon with the aim of different and specific applications: short-term scheduling which focuses to plan for less than a week like (De Souza Filho et al., 2013) and (Zhang et al., 2017), medium-term scheduling with the aim to plan for one to two weeks such as (Liao et al., 2018a), (Liao et al., 2018b) and (Chen et al., 2017), and long-term scheduling for horizons up to a month like (Meira et al., 2018), (Cafaro & Cerda, 2017), and (Cafaro et al., 2015).

**2.5. Revision in scheduling**

The main target of the proposed method is to achieve a model that can take changes in the operational environment into account and produces a plan that can adapt itself when events and incidents occur at the middle of operation.

Only a few number of researches tried to deal with events. Almost the first attempt to revise a scheduling is presented in (Relvas et al., 2007). In this research, alongside with an MILP model for a one-to-one pipeline, a mechanism is introduced to enable revision of the plan. This mechanism works in a simple way. It states that any time during the execution of the plan, if it is determined that there is a need to revise the initial plan, then the execution is stopped, changes are made in the parameters and a new plan will be generated.

Some other researches such as (Su et al., 2018) and (Mostafaei et al., 2016) have done focusing on the uncertainty of the model parameters. Another approach that has been investigated in (Taherkhani et al., 2017) in order to respond to the random nature of requirements is modeling with random variables and expectation calculations.

## 3. Proposed Model

In this section, we will describe the proposed mathematical model for pipeline planning, with the ability to be aware of and adapt to upcoming events. Our proposed model is based on the model presented in (Cafaro & Cerda, 2010). There are six main sets in this model: The set of pumping runs ($K$), injected batches ($I$), oil products ($P$), input terminals ($S$), and receiving destinations ($J$) are also included in the basic models. The last set ($E$) contains upcoming events in the planning horizon, ordered by time of occurrence. This means that the event $e$ must happen before $e'$ if $e < e'$. The Set $E$ has at least two elements $e_0$ and $e_{max}$, representing the start and the end of the planning horizon respectively.

*Table 1. Nomenclature*

| (a) Sets | |
|---|---|
| $K$ | chronologically ordered blocks of pumping runs |
| $I$ | product batches $(I = I^{\text{old}} \cup I^{\text{new}})$ |
| $I^{new}$ | new batches to be injected during the planning horizon |
| $I^{old}$ | old batches in the initial state of pipeline |
| $P$ | refined oil products |
| $P_i$ | refined oil product contained in old batch $i \in I^{old}$ |
| $S$ | pipeline input terminals |
| $J$ | pipeline output terminals |
| $J_p$ | distribution depots demanding product $p$ |
| $J_{p,s}$ | distribution depots demanding product $p$ to be supplied by an upstream source $s$ |
| $E$ | chronologically ordered events. It consists of at least two members, $e_0$ and $e_{max}$ representing the start and the end of planning |
| **(b) Parameters** | |
| $h_{max}$ | planning horizon length |
| $PV$ | total pipeline volume from the origin to the farthest depot |
| $T_e$ | occurrence time of event $e$. $\left(T_{e_0} = 0, T_{e_{max}} = h_{max}\right)$ |
| $cb_{p,j}$ | unit backorder penalty cost for tardy orders of product $p$ at output terminal $j$ |
| $cif_{p,p'}$ | total reprocessing cost of interface material involving products $p$ and $p'$ |
| $cin_{p,s}$ | pumping cost per unit of product $p$ at input terminal $s$ |
| $D_{min}, D_{max}$ | minimum/maximum delivery size to output terminals |
| $DL_{p,j}$ | minimum request of product $p$ at the output terminal $j$ |
| $DU_{p,j}$ | maximum amount of product $p$ that can be delivered to output terminal $j$ |
| $Q_{min}, Q_{max}$ | minimum/maximum batch injection size |
| $Q_{max,p}$ | maximum batch injection size for product $p$ |
| $SL_{p,s}$ | lowest amount of product $p$ to be shipped from input terminal $s$ |
| $SU_{p,s}$ | maximum amount of product $p$ that can be shipped from input terminal $s$ |
| $vb_{min,s}, vb_{max,s}$ | minimum/maximum product flow-rates at source $s$ |
| $W_i^0$ | initial volume of old batch $i$ |
| $\sigma_j$ | volume coordinate of output terminal $j$ from the pipeline origin |
| $\tau_s$ | volume coordinate of input terminal $s$ from the pipeline origin |
| **(c) Continuous variables** | |
| $BC$ | total backorder penalty cost for tardily meet product demands |
| $B_{p,j}$ | backorder of product $p$ at the output terminal $j$ |
| $C_k$ | completion time of run block $k$ (measured in time units) |
| $L_k$ | length of run block $k$ |
| $L_{k,s}$ | length of pumping run $k$ taking place at source $s$ |
| $D_{i,j}^{(k)}$ | volume of batch $i$ diverted to output terminal $j$ during run block $k$ |
| $DP_{i,j,p}^{(k)}$ | volume of product $p$ diverted from batch $i$ to terminal $j$ during run block $k$ |
| $F_{i,k}$ | upper coordinate of batch $i$ from the origin at time $C_k$ |
| $PC_k$ | total pumping cost during block $k$ |
| $Q_{i,s}^{(k)}$ | size of batch $i$ shipped from input terminal $s$ during run $k$ |
| $QP_{i,s,p}^{(k)}$ | size of batch $i$ containing product $p$ inputted from input terminal $s$ through run $k$ |
| $RC_i$ | total reprocessing cost for the interface between batch $i$ and the following $(i+1)$ |
| $W_{i,k}$ | size of batch $i$ at the end of run $k$ |
| **(d) Binary variables** | |
| $v_{i,s}^{(k)}$ | variable denoting that a portion of batch $i$ is injected from source $s$ through block $k$ |
| $x_{i,j}^{(k)}$ | variable denoting that a portion of batch $i$ is diverted to the depot $j$ during block $k$ |
| $y_{i,p}$ | variable denoting that batch $i$ contains product $p$ |
| $b_{k,e}$ | variable denoting that pumping run $k$ is executed completely inside of the interval $[T_{e-1}, T_e]$ |

### 3.1. Sequencing pumping run

A block $k \in K$ must start after the completion of the previous block $(k-1)$. Let $C_k$ denote the completion time and $L_k$ the duration of block $k$.

$$C_k - L_k \geq C_{k-1} \quad \forall k \in K (k > 1) \tag{1}$$

All pumping blocks should also be finished before the planning horizon $h_{max}$.

$$C_k \leq h_{max} \quad \forall k \in K \tag{2}$$

### 3.2. Aligning pumping runs with events

There is an assumption in the model that during the execution of a pumping run $k$, the environmental conditions and parameters remain stable. On the other hand, almost all events will change the parameters of the environment when they happen. Therefore, it should be prevented that an event occurs during the execution of a pumping run.

For this purpose, we define the binary variable $b_{k,e}$ in such a way that its value equals to 1, iff pumping run $k$ is executed completely between two consecutive events $e$ and $e-1$. It means that pumping run $k$ must be started after $T_{e-1}$ (the occurrence time of the event $e-1$) and must be finished before $T_e$ (the occurrence time of the event $e$). Figure 1 shows an example of this situation. In all other conditions, $b_{k,e}$ will be 0.

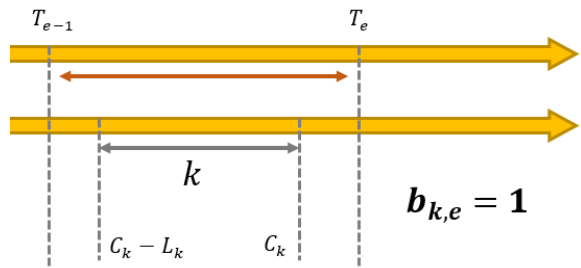

*Figure 1. An example of a pumping run fully executed between two consecutive events*

The following equations are added to the model to ensure the mentioned alignment. First, each pumping run $k$ should be executed no more than once. Then:

$$\sum_{e \in E} b_{k,e} = 1 \forall k \in K, e \in E \tag{3}$$

Second, if a pumping run is executed between two events, the start and finish time of the run should be placed between occurrence time of those events.

$$T_{e-1} b_{k,e} \leq C_k - L_k \forall k \in K, e \in E \tag{4}$$

$$C_k \leq T_e b_{k,e} + (1 - b_{k,e}) M_T \forall k \in K, e \in E \tag{5}$$

The value of constant $M_T$ should be larger than all time values in the model. So it can be chosen as $h_{max}$ or higher.

### 3.3. Injecting batches in different pumping runs

A pumping operation is identified by its corresponding pumping run $k \in K$ and source terminal $s \in S$. We define the binary variable $v_{i,s}^{(k)}$ to state that the batch $i \in I^{new}$ in injected to the pipeline from input terminal $s$ at the pumping run $k$. If $v_{i,s}^{(k)} = 1$, then the whole or at least some amount of batch $i$ is injected from source terminal $s$ during pumping run $k$. To ensure that no more than one batch is injected in the pipeline from each source in each pumping run, we have:

$$\sum_{i \in I} v_{i,s}^{(k)} \leq 1 \forall k \in K, s \in S \tag{6}$$

If for a certain pumping run $k$, the term $\left(\sum_{s \in S} \sum_{i \in I} v_{i,s}^{(k)}\right)$ is zero, it means that no batches are injected in that run. Since the size of set $(K)$ is normally chosen larger that the number of pumping runs in the optimal solution, this situation is probable for the last members of $(K)$ (named dummy runs).

The following Equation ensures that if a pumping run $k$ is dummy, all later runs also turn to dummy.

$$\sum_{s \in S} \sum_{i \in I} v_{i,s}^{(k)} \leq M \left(\sum_{s' \in S} \sum_{i' \in I} v_{i',s'}^{(k-1)}\right) \forall k \in K (k > 1) \tag{7}$$

The constant $M$ is a very large number, greater than all values that can be generated in the model. We will use it more in later equations in the model.

### 3.4. Batch injection and size relation

Every Batch $i$ is injected for the first time in a pumping run $k$ at a source terminal $s$. It is also possible that some portion of $i$ is appended to it in later pumping run $k'$ at another (or the same) terminal $s'$. We say that batch $i$ is injected in both runs $(k, s)$ and $(k', s')$.

Let $Q_{i,s}^{(k)}$ denote the size of injected product to batch $i$ in run $(k, s)$. $Q_{i,s}^{(k)}$ will be positive iff there is actually a real injection in batch $i$ at run $(k, s)$. This means that $v_{i,s}^{(k)}$ and $Q_{i,s}^{(k)}$ will be either zero or positive simultaneously:

$$Q_{min} v_{i,s}^{(k)} \leq Q_{i,s}^{(k)} \leq Q_{max} v_{i,s}^{(k)} \forall i \in I, s \in S, k \in K \tag{8}$$

Parameters $Q_{min}$ and $Q_{max}$ present the global lower and upper bound for the volume of product that can be injected during a pumping run.

### 3.5. Choosing lengths of pumping runs

Let $L_{k,s}$ be the length of the run $(k,s)$, means the time length during run $k$ that terminal $s$ is actually pumping products. The whole run $k$ can be finished only if all of its sub-runs $(k,s)$ is finished before. So the length of run $k$ will be given by $L_k = \max_{s \in S}(L_{k,s})$ or:

$$L_k \geq L_{k,s} \forall k \in K, s \in S \tag{9}$$

The injected volume of refined product is also related to the length of pumping run and the rate at which a terminal can perform injection. The total volume of product injected to the pipeline during run $(k,s)$ is equal to $\sum_{i \in I} Q_{i,s}^{(k)}$. In the base model from (Cafaro & Cerda, 2010) the following equation is presented, where $vb_{min,s}$ and $vb_{max,s}$ stand for the minimum and maximum pumping rate of terminal $s$.

$$vb_{min,s} L_{k,s} \leq \sum_{i \in I} Q_{i,s}^{(k)} \leq vb_{max,s} L_{k,s} \forall k \in K, s \in S \tag{10}$$

In this model, we want to take the effects of events into consideration. The next section discusses the problem and presents the proposed technique to address this requirement.

### 3.6. Choosing parameters based on event occurrences

When an event occurs, it may affect the conditions of pipeline and its surroundings. In our model, this will be done by changing the value of parameters. Therefore, the value of a parameter (same as $vb_{min,s}$ and $vb_{max,s}$ from Equation 8) should be chosen according to event occurrences. Since the parameters are assumed to be fixed when solving an MILP model, our technique to address this need is to have a vector of all values that a parameter can take over the planning horizon and using the respective value in the constraint expressions.

For example, suppose there is a predetermined maintenance program at terminal $s$ (represented by event $e$) that will cause the maximum pumping rate at the terminal to decrease to half of the previous value, and it remains the same until the maintenance end (represented by event $e'$). In such a case, we consider three values for the parameters $vb_{max,s}$: the value before the occurrence of $e$, the value between $e$ and $e'$, and the value after $e'$.

Let us denote the value of an arbitrary parameter $A$ between events $e-1$ and $e$ by the notation $A^{\langle e \rangle}$. So the three parameters mentioned before will be denoted by $vb_{max,s}^{\langle e \rangle}$, $vb_{max,s}^{\langle e' \rangle}$, and $vb_{max,s}^{\langle e_{max} \rangle}$.

According to these assumptions and notations, we rewrite the Equation 10:

$$vb_{min,s}^{<e>} L_{k,s} - M(1 - b_{k,e}) \leq \sum_{i \in I} Q_{i,s}^{(k)} \leq vb_{max,s}^{<e>} L_{k,s} + M(1 - b_{k,e}) \forall k \in K, s \in S, e \in E(e > 0) \tag{11}$$

The term $M(1 - b_{k,e})$ is called *event guard* and will be noted as $EG_{k,e}$ hereafter. It is used to ensure that if the pumping run $k$ is executed between events $e-1$ and $e$ (means $b_{k,e}$), only the corresponding equation (containing parameters $A^{\langle e \rangle}$) takes effect and the other ones go trivial. That will be held because an infinity term is added to the larger side of the inequality.

### 3.7. Tracking the batch size over time

Let $W_{i,k}$ be the size of a batch $i$ at the time $C_k$. If the volume of batch $i$ is changed during run $k$, it must be either an injection from source $s$ or a delivery to depot $j \in J$.

$$W_{i,k} = W_{i,k-1} + \sum_{s \in S} Q_{i,s}^{(k)} - \sum_{j \in J} D_{i,j}^{(k)} \forall i \in I, k \in K \tag{12}$$

Only the injections that are performed at upstream sources can cause older batches to move forward. So the amount of product that can be delivered to depot $j$ is at most equal to the amount of injections performed at upstream sources $(\tau_s < \sigma_j)$.

$$D_{i,j}^{(k)} \leq \sum_{s \in S: \tau_s < \sigma_j} \sum_{i' \in I} Q_{i',s}^{(k)} \forall i \in I, j \in J, k \in K \tag{13}$$

### 3.8. Tracking the batch location over time

Let $F_{i,k}$ denote the upper volumetric coordinate of batch $i$ at the time $C_k$. So $F_{i+1,k}$ also represents the upper coordinate of batch $(i+1)$, which is equal to the lower coordinate of batch $i$ at the same time. Since the difference between lower and upper coordinates of a batch $i$ at the end of run $k$ is equal to its size $W_{i,k}$, we have:

$$F_{i,k} - W_{i,k} = F_{i+1,k} \forall i \in I, k \in K \tag{14}$$

To prevent backward movement of batches during the time:

$$F_{i,k-1} \leq F_{i,k} \forall i \in I, k \in K \tag{15}$$

Finally to ensure that all batches stands into pipeline and no batches will ever be located outside the pipeline coordination, Let $PV$ stands for the total pipeline volume, then the upper coordinate of any batch $i$ must never be greater than $PV$ and its lower coordinate must be non-negative.

$$F_{i,k} \leq PV \forall i \in I, k \in K \tag{16}$$

$$F_{i,k} - W_{i,k} \geq 0 \forall i \in I, k \in K \tag{17}$$

### 3.9. Preventing empty spaces and overflows

The pipeline must remain full of products at any time over the planning horizon. Therefore, the sum of volumes of all batches in the pipeline must be equal to $PV$:

$$\sum_{i \in I} W_{i,k} = PV \forall k \in K \tag{18}$$

Also the amount of injected and delivered products must be equal at the end of any run $k$:

$$\sum_{i \in I} \sum_{s \in S} Q_{i,s}^{(k)} = \sum_{i \in I} \sum_{j \in J} D_{i,j}^{(k)} \forall k \in K \tag{19}$$

### 3.10. Injections and deliveries feasibility

Refined product can be injected into batch $i$ during pumping run $(k, s)$ only if the following two conditions are satisfied:

(a) Batch $i$ has reached the location of the source $s$ ($\tau_s$) before the run $k$ begins:

$$F_{i,k-1} \geq \tau_s v_{i,s}^{(k)} \quad \forall i \in I, s \in S, k \in K \tag{20}$$

(b) Also batch $i$ has not completely passed the location of source $s$ before the run $k$:

$$F_{i,k-1} - W_{i,k-1} \leq \tau_s + M(1 - v_{i,s}^{(k)}) \quad \forall i \in I, s \in S, k \in K \tag{21}$$

### 3.11. Delivery from batches to depots

Same as the previous section, delivery of product from batch $i$ to depot $j$ is also possible, only if:

(a) The batch $i$ reaches the location of depot $j$ (denoted by $\sigma_j$) before the end of run $k$, i.e. $F_{i,k} \geq \sigma_j$

(b) The batch $i$ has not passed the location of depot $j$ before the run $k$ begins. It means that the lower coordination of the batch must not exceed the location of the depot, i.e. $F_{i,k-1} - W_{i,k-1} \leq \sigma_j$

To fulfill these conditions, let $x_{i,j}^{(k)}$ be a binary variable denoting that some product is delivered from batch $i$ to depot $j$ during run $k$. When $x_{i,j}^{(k)} = 0$, nothing should be transferred from batch $i$ to depot $j$. It means $D_{i,j}^{(k)} = 0$.

$$D_{min,j}^{<e>} x_{i,j}^{(k)} \leq D_{i,j}^{(k)} \leq D_{max,j}^{<e>} x_{i,j}^{(k)} \quad \forall i \in I, j \in J, k \in K, e \in E(e > 0) \tag{22}$$

where $D_{min,j}^{\langle e \rangle}$ and $D_{max,j}^{\langle e \rangle}$ are bounds on the amount of product that can be delivered to depot $j$ in timespan $[T_{e-1}, T_e]$.

The following two constrains are used to meet the two conditions mentioned before in this section, respectively.

$$F_{i,k} \geq \sigma_j x_{i,j}^{(k)} \quad \forall i \in I, j \in J, k \in K \tag{23}$$

$$F_{i,k-1} - W_{i,k-1} \leq \sigma_j + M(1 - x_{i,j}^{(k)}) \quad \forall i \in I, j \in J, k \in K \tag{24}$$

As mentioned before, the prerequisite for delivering product from batch $i$ to depot $j$ at run $k$ is that the following inequalities are true:

$$F_{i,k-1} - W_{i,k-1} < \sigma_j \leq F_{i,k} \tag{25}$$

In this case, an upper bound on the volume of product that can be delivered is $\sigma_j - (F_{i,k-1} - W_{i,k-1})$, unless more product is injected in batch $i$ from upstream sources at the same time that delivery is in action.

From the other hand, delivery of product can be done at any depot that is in touch with the batch during the run. So the total amount of deliveries from batch $i$ to all the depots from the origin to depot $j$, could not be more than the initial available volume of batch $i$ $\left(\sigma_j - (F_{i,k-1} - W_{i,k-1})\right)$ added by the total amount of injections from upstream sources to batch $i$. Then:

$$\sum_{j'=1}^{j} D_{i,j'}^{(k)} \leq \sigma_j - (F_{i,k-1} - W_{i,k-1}) + \sum_{\substack{s \in S \\ \tau_s < \sigma_j}} Q_{i,s}^{(k)} + (PV - \sigma_j)(1 - x_{i,j}^{(k)}) \quad (26)$$
$$\forall i \in I, j \in J, k \in K$$

### 3.12. Assigning products to batches

At most, only one product can be assigned to batch $i$. Let $y_{i,p}$ be a binary variable denoting that batch $i$ contains product $p$ whenever $y_{i,p} = 1$. Then:

$$\sum_{p \in P} y_{i,p} \leq 1 \quad \forall i \in I \quad (27)$$

If no products is assigned to batch $i$, then we say that $i$ is fictitious and will not be injected in the pipeline in real life. It means that there is no pumping run $(k,s)$ injecting any product to batch $i$. Then we have $v_{i,s}^{(k)} = 0$ for all pumping runs.

In the other hand, if there is a product assignment to batch $i$, it means that this batch is real and at least one pumping run must inject some products to it. This condition can be mathematically written as follows,

$$\sum_{p \in P} y_{i,p} \leq \sum_{s \in S} \sum_{k \in K} v_{i,s}^{(k)} \leq M \sum_{p \in P} y_{i,p} \quad \forall i \in I^{new} \quad (28)$$

Some products should not be injected next to each other because of contamination. Let us assume that $(p, p')$ stands for a forbidden sequence of products. Then,

$$y_{i-1,p} + y_{i,p'} \leq 1 \quad \forall i \in I^{new}(i > 1) \quad (29)$$

### 3.13. Product injection

The whole amount of products injected into batch $i$ in all pumping runs must be from a single type of product $p$, such that $y_{i,p} = 1$.

$$QP_{i,s,p}^{(k)} \leq Q_{max} y_{i,p} \quad \forall i \in I, s \in S, k \in K, p \in P \quad (30)$$

$$\sum_{p \in P} QP_{i,s,p}^{(k)} = Q_{i,s}^{(k)} \quad \forall i \in I, s \in S, k \in K \quad (31)$$

### 3.14. Product Delivery

Same as the previous section, The whole amount of delivered material from batch $i$ to different depots must be from a single type of product $p$ such that $y_{i,p} = 1$.

$$DP_{i,j,p}^{(k)} \leq D_{max} y_{i,p} \quad \forall i \in I, j \in J, k \in K, p \in P \quad (32)$$

$$\sum_{p \in P} DP_{i,j,p}^{(k)} = D_{i,j}^{(k)} \quad \forall i \in I, j \in J, k \in K \quad (33)$$

### 3.15. Fulfilling product demands

The model should provide the demanded amount of products for all depots. Let say that the total amount of product $p$ needed at the depot $j$ is $DL_{p,j}$. To let the model continue to work and be feasible, even if it does not fulfill the demand $DL_{p,j}$, we introduce variable $B_{p,j}$ as the amount of product shortage. It will be used in the cost function later.

At last, to consider the upper bound of storage capacity for product $p$ at depot $j$, we use parameter $DU_{p,j}$. So the total amount of product $p$ delivered to depot $j$ must be lower that this bound.

$$DL_{p,j} - B_{p,j} \leq \sum_{k \in K} \sum_{i \in I} DP_{i,j,p}^{(k)} \leq DU_{p,j} \forall p \in P, j \in J \tag{34}$$

### 3.16. Initial state of problem

To initialize the starting position and order of batches in the pipeline, let $W_i^0$ be the volume of the old batch $i \in I^{old}$. Then, the upper coordinate of batch $i \in I^{old}$ is equal to the sum of initial volumes for all succeeding batches $i' \in I^{old}$ ($i' > i$) plus its own volume:

$$F_{i,0} = \sum_{i' \in I^{old}(i' \geq i)} W_{i'}^0 \forall i \in I^{old} \tag{35}$$

Moreover, it is also known the product $P_i$ contained in every old batch $i \in I^{old}$

$$y_{i,P_i} = 1 \forall i \in I^{old} \tag{36}$$

### 3.17. Objective function

The main goal of the proposed method is to increase the adaptability of the model to changes in the operating environment. Therefore, in this method, in addition to trying to achieve a program that will deliver the demand of depots with at the lowest cost of pumping ($PC$) and reprocessing of interface material ($RC$), the cost of shortage in the needs for the products ($BC$) must be reduced as much as possible.

$$minz = \sum_{i \in I} RC_i + \sum_{k \in K} PC_k + BC \tag{37}$$

The following three equations introduce the terms of the objective function:

$$RC_i \geq \text{cif}_{p,p'}(y_{i,p} + y_{i+1,p'} - 1) \forall i \in I; p, p' \in P \tag{38}$$

The variable $RC_i$ stands for the total cost of reprocessing interface material between batches $i$ and $i+1$ and the parameter $cif_{p,p'}$ is the cost for reprocessing the interface material resulting from blending of products $p$ and $p'$.

$$PC_k = \sum_{i \in I} \sum_{s \in S} \sum_{p \in P} \text{cin}_{p,s} QP_{i,s,p}^{(k)} \forall k \in K \tag{39}$$

The variable $PC_k$ is the total pumping cost at run $k$ and the coefficient $\text{cin}_{p,s}$ stands for the cost of pumping a unit volume of product $p$ in the source $s$.

$$BC = \sum_{p \in P} \sum_{j \in J} \text{cb}_{p,j} B_{p,j} \tag{40}$$

The parameter $\text{cb}_{p,j}$ stands for the cost of failing to provide a single unit of product $p$ at depot $j$ on schedule and $B_{p,j}$ (as introduced in Equation 34) is the amount of backorder of product $p$ at depot $j$.

## 4. Results & Discussion

In this section, the results of the implementation of the proposed model and its performance comparison with previous methods are reported. The study case considered for this comparison is based on an example mentioned in (Cafaro & Cerda, 2010), which has been modified to meet the requirements of this comparison.

### 4.1. Example configuration

In this example, the operational planning of a pipeline network is carried out, which includes two source stations $S1$ and $S2$ and three depots $D1$, $D2$ and $D3$. This pipeline transfers three refined products $A$, $B$ and $C$ from the sources to depots. Stations $S1$ and $S2$ are located at volume coordinates 0 and 40 from the beginning of the pipeline, respectively. Stations $D1$, $D2$, and $D3$ are also located at volume coordinates 20, 60, and 80 from the pipeline origin. Between each two consecutive stations, there is a pipeline with a volume of 20 units.

The available amount of the products at the sources and the demands of the depots for the products are shown in Table 2. Also, in Table 3, the cost of recycling the interface between each pair of products can be seen.

*Table 2. Product supplies and demands*

|   | Supplies | | Demands | | |
|---|---|---|---|---|---|
|   | S1 | S2 | D1 | D2 | D3 |
| A | 50 | 20 | 60 | 60 | - |
| B | 80 | 60 | - | - | 100 |
| C | 30 | 80 | - | 60 | - |

*Table 3. Interface costs (in 100$)*

| Predecessor | Successor | | |
|---|---|---|---|
|   | A | B | C |
| A | - | 22 | 23 |
| B | 24 | - | 21 |
| C | 30 | 32 | - |

Table 4 contains the cost of pumping one volume unit of each product in each of the sources. The feasible range for pumping rate for all products in all sources is set to be 1.2 (units/h).

*Table 4. Pumping cost in each source (in 100$ per volume unit)*

| Product | Source | |
|---|---|---|
|   | S1 | S2 |
| A | 29 | 14.5 |
| B | 34 | 17.0 |
| C | 24.5 | 49.0 |

The initial state of the pipeline includes four batches $B1$, $B2$, $B4$ and $B5$. The product and the initial volume of these batches is shown in the top row of Figure 2.

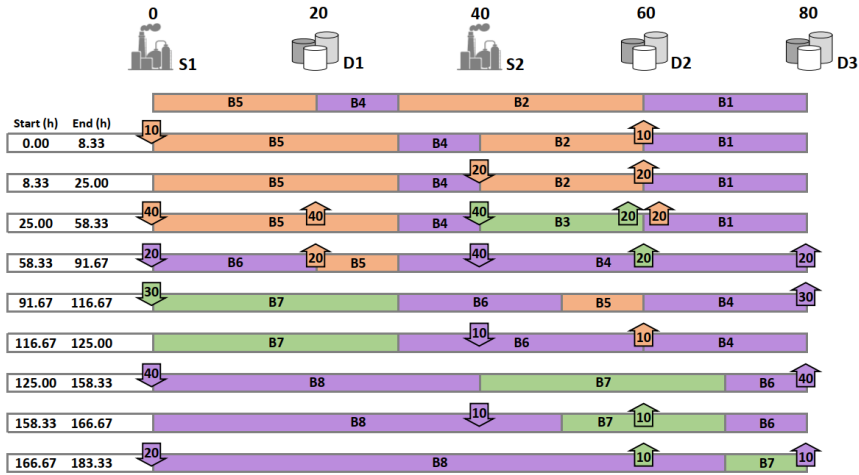

Figure 2. Optimal pipeline schedule for the example when no events occur

### 4.2. Optimal solution in absence of incidents

When no events occur during the planning horizon, the optimal schedule for transporting products and meeting the demands of the depots is as shown in Figure 2. In this case, we will have 9 pumping runs $k1 - k9$. A distinct row is drawn in Figure 2 for each of pumping runs. The start time, completion time and duration of each run can be seen at the left of each row and the volume of the product injected into the pipeline or received from it during the run is also drawn as arrows. The color pattern of each batch indicates the product related to it, as is defined in the legend of the figure.

As can be seen, the time required to execute the optimal schedule is 183.33 hours. The plan is as follow:

In the first run ($k1$), which starts from time 0 and continues for 8.33 hours, 10 units of product $A$ are injected from station $S1$ to batch $B5$. As the batches move forward in the pipeline, depot $D2$ takes 10 units of product $A$ from batch $B2$. In the second run ($k2$) which takes 16.67 hours long, station $S2$ performs the pumping operation and adds 20 units of product $A$ to batch $B2$ and exactly the same amount of product from the same batch is received again in depot $D2$. Run $k3$ takes 33.33 hours and during this run, 40 units of product $A$ are injected to batch $B5$ from station $S1$ and the same amount is received in depot $D1$. As there is no pumping interaction and batch overlap during this run, station $S2$ can also pump 40 units of product $C$ in the form of a new batch $B3$ simultaneously. This injection causes the depot $D2$ to fully receive the remaining 20 units of product $A$ in $B2$ and, in addition, withdraw 20 units of product $C$ from batch $B3$. Pumping Run $k4$, which ends at 91.67 hours, includes two injections of product $B$: 20 units from station $S1$ to the new batch $B6$ and 40 units from station $S2$ to $B4$ that had reached the station at the end of the previous run. As a result of these two injections, 20 units of the product $A$ from $B5$ enter depot $D1$, the remaining 20 units of the product $C$ in the batch $B3$ are received in depot $D2$, and 20 units of product $B$ are added to the inventory of depot $D3$.

As seen in Figure 2, during runs $k5$ to $k8$, only one injection and one delivery take place in each run. In the last pumping run ($k9$), station $S1$ injected 20 units of product $B$ into batch $B8$, and as a result, depots $D2$ and $D3$ received respectively, 10 units of product $C$ from batch $B7$ and 10 units of product $B$ from batch $B6$. The final state of the pipeline, 183.33 hours after the start of the program is as follow: the batch $B8$ containing product $B$ is location from the beginning of the pipeline to the volumetric coordinate 70, and batch $B7$ is placed from the coordinates 70 to 80, containing 10 units of the product $C$.

### 4.3. Adding incidents to the example

It is mentioned earlier that incidents usually lead to unwanted changes in the schedule and as a result, the time and cost of the plan increase according to the severity of the incident. In this example, we consider a case where a maintenance operation will be executed on the pipeline between $S1$ and $D1$ after 100 hours from the start of the planning horizon. Therefore, during this operation, it will not be possible to inject any product from the source $S1$. We also assume that the maintenance will be finished in 30 hours.

In none of the previous pipeline planning methods, incidents and changes in the operational environment have been considered. So there is no way to add events to the old models. When using these models in the real world, when an incident occurs, the only way is to stop the plan execution, update the numbers and parameters of the model according to the new state of the environment, and run the model again. If we're lucky, the model will probably lead to a new plan that will definitely cost more than before. But in many cases, it is not possible to achieve any plan due to the constraints of the model, and practically, the model that worked well before the occurrence of the event, after that only reports infeasible results.

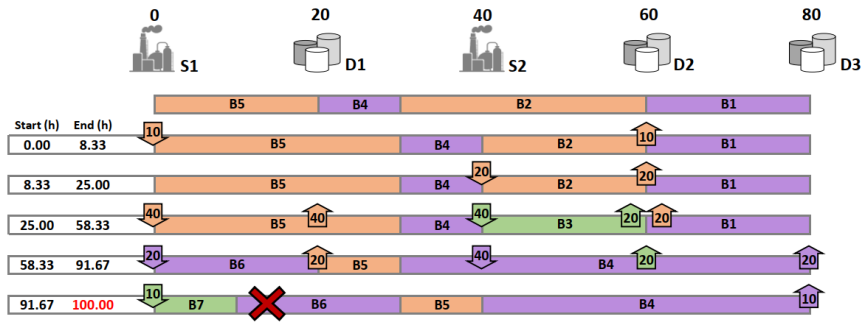

Figure 3. State of operational environment at time $t=100h$

However, in this example, we proceed through the mentioned method. From time $t = 0$ to $t = 100h$, it will proceed according to the previous plan and the batches will be injected in exactly the same amount and in the same order. Therefore, at $t = 100h$ and at the moment of starting the maintenance operation, the situation is as shown in Figure 3. After that, due to the reduction of maximum rate for injection in the pipeline to zero, there is no way to inject any new batches in the pipeline and practically, the state of the pipeline remains constant for next 30 hours (which clearly shows the waste of time and resources). From $t = 130h$, the planned program continues and the rest of the steps and pumping runs are done according to what was assumed in plan.

Finally, at $t = 213.33h$, which is 30 hours later than the optimal schedule, the distribution of the products finishes. Step by step execution of this plan has been drawn in Figure 4.

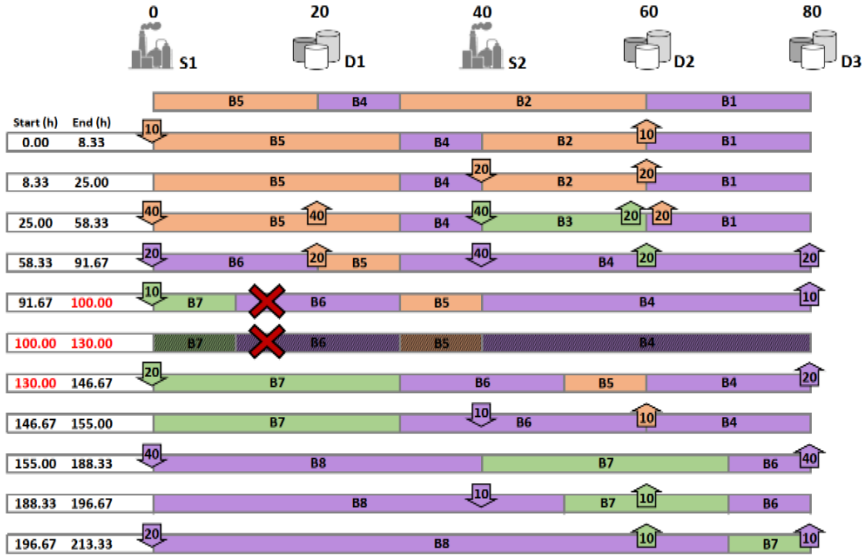

Figure 4. The output plan of older methods when the incident happens

### 4.4. Performance of the proposed method in facing the incident

The advantage of the proposed method is that during planning, it pays attention to the future changes of the environment and can consider different values for the parameters. It is also possible to consider predictable events (such as the maintenance program mentioned in this example) from the beginning.

In the proposed method, we consider two changes compared to previous methods during modeling. First, we define the set of events as follow:

$$E = \{e_0, e_1 \text{ (Maintenance starts)}, e_2 \text{ (Maintenance ends)}, e_{max}\} \quad (41)$$

As explained earlier, in this set, two events $e_0$ and $e_{max}$ are included to represent the beginning and end of the time horizon. Events $e_1$ and $e_2$ also represent the start and end of the maintenance program. Therefore, the values of the parameters $T_{e_1}$ and $T_{e_2}$ are considered to be $100h$ and $130h$, respectively.

The second change is to define a vector of values for parameters $vb_{min,s_1}$ and $vb_{max,s_1}$. The values in this vector are as follow (values are in units per hour):

$$\begin{aligned} vb_{min,s_1}^{\langle e_1 \rangle} &= 1.0 & vb_{min,s_1}^{\langle e_2 \rangle} &= 0.0 & vb_{min,s_1}^{\langle e_{max} \rangle} &= 1.0 \\ vb_{max,s_1}^{\langle e_1 \rangle} &= 1.2 & vb_{max,s_1}^{\langle e_2 \rangle} &= 0.0 & vb_{max,s_1}^{\langle e_{max} \rangle} &= 1.2 \end{aligned} \quad (42)$$

In this way, the model realizes that in the interval $[T_{e_1}, T_{e_2}]$, no pumping should be performed from station $S1$. For other parameters, a vector of length $|E| - 1$ should be used too, but all the values in them will be the same because the only parameters changing through time in this example are $vb_{min,s_1}$ and $vb_{max,s_1}$

The result of our proposed model has a much better performance than the old methods. As can be seen in Figure 5, this method already looks into the future and plans for the time that station $S1$ is not pumping, in such a way that the distribution of products from other stations continues.

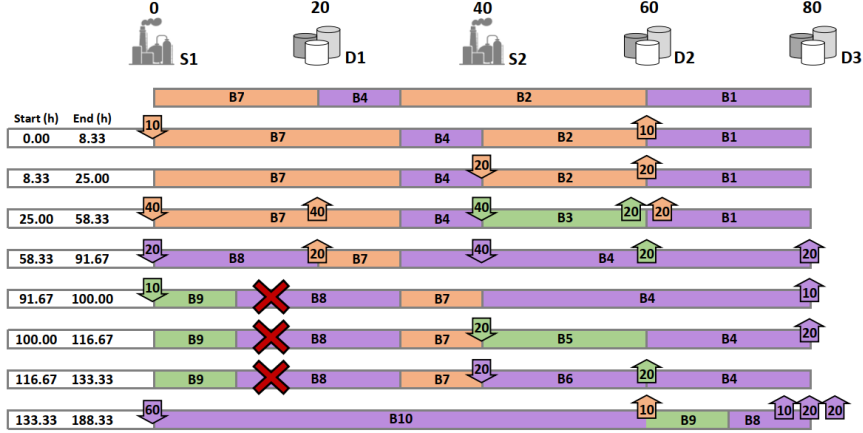

*Figure 5. The output plan of our proposed model*

The output of our method is same as the output of the previous methods until the moment $t = 100h$. It pumps batches of different products in the pipeline in the same order and with the same amount. But the divergence starts from time $T_{e_1}$. Considering that from this moment until $T_{e_2}$ the pumping rate at station $S1$ will be zero, the first part of the pipeline (from the beginning to the volume coordinates of 40 units) remains constant. But in the same period of time, station $S2$ starts pumping products that satisfy part of the demands of depots $D2$ and $D3$. In this way, when in $T_{e_2}$ i.e. 30 hours later, the station $S1$ starts pumping again, the remaining amount of products needed in $D2$ and $D3$ is reduced compared to the start of the maintenance operation. So by pumping next batches, it can be reached the desired condition faster and therefore with lower cost.

More precisely, the output of our model consists of 8 pumping run ($k1 - k8$), where the runs $k1$ to $k4$ are performed exactly the same as in the normal case. The only difference is that batch $B5$ in normal case is numbered as $B7$ here, because of two intermediate batches that will be injected to the middle of the pipeline in future runs.

The fifth run ($k5$), which is forced to be left unfinished in the previous methods, ends consciously at time $T_{e_1}$ in this model. Run $k6$ begins exactly when event $e_1$ occurs and thus, station $S1$ performs no pumping at this run. During the 16.67 hours of run $k6$, station $S2$ pumps 20 units of product $C$ as a new batch called $B5$ and on the other hand, depot $D3$ withdraws 20 units of product $B$ from batch $B4$. In run $k7$, station $S1$ is still inactive so the station $S2$ again injects a new batch named $B6$ containing 20 units of product $B$ into the pipeline and batch $B5$ is completely delivered to the depot $D2$.

Run $k8$ starts after the maintenance operation is finished, and therefore, station $S1$ is reactivated in this run. Station $S1$ starts pumping a large new batch $B10$ with a volume of 60 units of product $B$. In this way, the depot $D2$ takes the last part of its demand, i.e. 10 units of the product $A$ from the batch $B7$, and the depot $D3$ takes three consecutive withdrawals. The remaining volume of the batches $B4$ and $B6$ plus a portion of batch $B8$ -which are equal to 50 units of product $B$ in overall- are delivered to the depot $D3$. The final state of the pipeline is very similar to the final state in the normal case. There are 70 units of product $B$ and 10 units of product $C$ in the pipeline.

### 4.5. Conclusion

The execution time of the output plan of our proposed model is surprisingly equal to the required time in normal case. It means that sometimes even severe incidents like pipeline shutdowns, can be managed easily and with the lowest overhead cost.

Our proposed model generates a plan which is 30 hours shorter than the time required to execute the program in the previous methods, due to the fact that the pipeline is not stopped totally during the maintenance operation. The total demand of all depots is fulfilled in only 183.33 hours which is 16\% faster than the output plan of older methods, despite the possibility of older methods to go totally infeasible in many other cases.